# Reduced-Order Modeling of Bolt Loosening: Application to a Pair of Oscillators Under Transverse Shock Excitation

Qirui He[1], Rui Wang[1,2], Matthew J. Alexander[1], Keegan J. Moore[1]*

[1]*Daniel Guggenheim School of Aerospace Engineering, Georgia Institute of Technology, Atlanta, GA 30332*
[2]*Department of Mechanical Engineering, Rice University, Houston, TX 77005*

*Corresponding author: (K.J. Moore)
E-mail address: kmoore@gatech.edu

## ABSTRACT

The safety and integrity of engineered structures are critically dependent on maintaining sufficient preload in their bolted joints. This preload can be dynamically lost due to sustained vibrations or sudden shock that are large enough to induce slip in the threads. While high-fidelity finite element simulations and analytical methods can accurately model the loss of preload for a single, their prohibitive computational expense and complexity render them unfeasible for analyzing large-scale structures with many bolts. This creates a critical need for reduced-order models that capture the essential physics of loosening while remaining computationally efficient. This paper introduces a reduced-order modeling methodology for predicting the loosening of bolted lap joints subjected to transverse shock excitation. The core idea is to treat the bolt tension as a dynamic degree-of-freedom that governs the effective properties of the joint through tension-dependent stiffness and damping that couple the components together. The methodology is applied to a pair of oscillators coupled by with a single lap joint with a strain-sensing bolt. Three different sets of experimental measurements are used to interrogate the dynamics of the system. Mathematical models are identified for the joint stiffness and damping and the instantaneous tension, which are combined with the equations of motion for the oscillators to simulate and reproduce the experimental measurements. Ultimately, the results validate the treatment of bolt tension as a dynamic degree-of-freedom, such that the methodology provides an effective framework for predicting loosening behavior in bolted joints.

**Keywords:** bolts, bolt loosening, reduced-order modeling, self-loosening, joint mechanics, bolted joints

## 1. Introduction

Bolted joints are the most common joining mechanism and are used to connect components in nearly every industry. Generally, structures contain hundreds if not thousands of bolted joints and loosening of even just one bolt can lead to catastrophic failure. For example, the loosening of a single screw caused an EQ-4B Global Hawk to lose stability and crash in 2011 [1]. In 2016, a Union Pacific train derailed in Mosier, Oregon due to a single loose bolt, which released 42,000 gallons of crude oil into the local environment [2–4]. Consequences of loose bolts are not just limited to extreme engineering failures. In fact, a survey of 103 playgrounds in the US Midwest found that 29.1% had loose fasteners and were dangerous for children [5–7]. Bolt loosening is generally split into non-rotational and rotational loss of preloads [8]. Non-rotational loosening focuses on the effects that do not involve relative motion between the threads of the joint [9–15] and is not considered further in this work. Rotational loosening occurs when the mating threads slide against each other, releasing and dissipating the elastic energy that was stored in the joint from tightening. Rotational loosening often leads to catastrophic failure do to a sudden and complete loss of preload, such that it is considered one of the worst failure mechanisms in many industries [8]. Throughout this paper, we refer to rotational loosening as just loosening.

Modeling of bolted joint loosening has become popular over the past three decades with analytical methods and high-fidelity finite element (FE) models representing the bulk of models. Researchers have developed various analytical models to understand and predict loosening in bolted joints. Hess and his students constructed reduced-order models for loosening under axial excitation where the bolt threads are modeled as free masses trapped between two inclines [16,17]. They extended their model to capture the dynamics of the full bolt and applied that model to evaluate the effect of the location of a bolted joint in a beam [20]. Further contributions come from Nassar and Yang [21–27], who formulated integral equations to describe the different friction forces and the bolt tension. Building on this, Nassar and Abboud [28] developed an improved stiffness model for bolted joints, while Yang et al incorporated the effects of plasticity in their models [29]. More recently, a new criteria for loosening was proposed by Fort et al. [30,31], which establishes a relationship between the transverse deflection of the bolt and the relative transverse displacement across joint.

The use of finite element models for studying loosening often presents significant challenges, particularly the need for highly refined meshes to accurately represent the thread geometry and contact mechanics. Consequently, many finite element models are limited to only a single bolt and are ran for only a few cycles of vibrations. One of the first FE model was developed by Pai and Hess [32], which was able to reproduce different loosening motions measured experimentally. Jiang et al. [33] later implemented both elastic and plastic effects into loosening, but neglected the helical geometry to simplify the problem. Subsequent models by Izumi et al. [11,34] and Zhang et al. [35] offered improvements, but their mesh densities were insufficient for detailed analysis of local contact conditions. A more detailed examination of local slip conditions within the threads and at the bolt head during self-loosening was conducted by Dinger and Friedrich [36]. The distinct effects of rotational and axial excitations were explored in separate studies in [37,38] using models with approximately 200,000 nodes. Furthermore, the influence of thread wear has been predicted by Zhang et al. [14,15] using models of increasing complexity, with 80,546 and 432,529 nodes, respectively. More recent investigations into anti-loosening performance have consistently employed models with around 500,000 elements [40–45]. In a particularly demanding simulation, Liu et al. [46] studied the impact of torsional excitation on self-loosening with an FE model comprising 1,385,142 elements. Although these and other studies have successfully created FE models that can replicate experimental findings, they are almost invariably restricted to the behavior of a single bolt. The intricate geometry of bolt threads makes them exceptionally difficult to mesh, demanding extremely small elements that can lead to convergence problems in standard FE software. When these models do converge, the vast number of degrees-of-freedom resulting from the fine mesh leads to computationally expensive and time-consuming simulations. Consequently, performing high-fidelity simulations of loosening across entire structures, which could contain numerous bolted joints, remains impractical with current computational capabilities.

The significant computational expense associated with current modeling techniques highlights the necessity for developing reduced-order models (ROMs). These models aim to capture the primary structural and dynamic consequences of joint loosening while requiring substantially less computational effort. In earlier work, Moore [47] presented a ROM framework for simulating the loosening behavior of a threaded joint connecting two axially aligned rods under shock loading [48–51]. This methodology models the stiffness of the joint as a function of its current torque, which is itself treated as a dynamic degree-of-freedom. A first-order, homogeneous ordinary differential equation was employed to describe the loss of torque over time. Concurrently, the relationship between the stiffness of the joint and its torque was represented by a function that increases monotonically. When this ROM was integrated into an FE model for the connected rods, it successfully replicated experimental results. The model accurately predicted both the strains recorded along the rods and the overall decrease in torque for various initial torque values and shock magnitudes. This ROM was later applied by Aldana and Moore to a computational system consisting of three axial rods connected with two axially threaded joints. The joints used the same models for joint stiffness, but the torques were treated as separate degrees of freedom. The response of the system was simulated for shock excitation with varying amplitude and initial starting torques. The results revealed that when the two joints were initially secured with identical torques, they both loosened completely when subjected to forces of sufficient magnitude. Furthermore, the research showed that the complete loosening of only a single joint occurred exclusively when the two joints had different initial torques.

Building on the success of our past work [47,52], the objective of this paper is to develop and apply a reduced-order modeling methodology for loosening of bolts in lap joints. The core assertion in our method is that the tension is a dynamic degree of freedom that governs the effective stiffness and damping of the joint. The system considered is a pair of harmonic oscillators connected by a single lap joint with a single bolt. Given the that the surrogate hypothesis for bolted joints has been shown to be true [53], the system studied here could be used as a surrogate for other bolted joints. The remainder of this paper is divided as follows. Section 2 introduces the proposed modeling approach, the experimental system, the measurement setups, and the bolt calibration. Section 3 investigates each measurement setup and leverages the measured data to form a mathematical model for the tension-dependent stiffness and damping of the joint and a differential equation governing the instantaneous tension. Section 4 combines all identified models together to form a final model, which is used to simulate and reproduce the experimental measurements for a full range of preloads. Finally, some concluding remarks are provided in Section 5.

## 2. Proposed Approach and Experimental System

### 2.1. Proposed Modeling Approach

We consider the loosening of bolted lap joints subjected to tranverse vibrations induced by impact excitation. The core of the modeling approach is to assert that the bolt tension, $T(t)$, is a dynamic degree-of-freedom (DOF) that governs the stiffness and damping properties of the joint. Our previous research [47] implemented a similar modeling approach for the torque. Although friction is present in bolted joints and the main source of dissipation, this research focuses on producing reduced-order models for the changes in stiffness and viscous damping that occurs with the loss of bolt preload. Reduced-order modeling of effects of friction, including both micro- and macro-slip, are left open for future research. To this end, we consider the system shown in Fig. 1 comprised of two elastic members coupled together using a bolt and nut. The joint is modeled using a spring with stiffness $k(T)$ and a viscous damper with damping coefficient $c(T)$, such that the elements are linear with respect to the relative motion of the structure but may nonlinearly depend on the instantaneous tension.

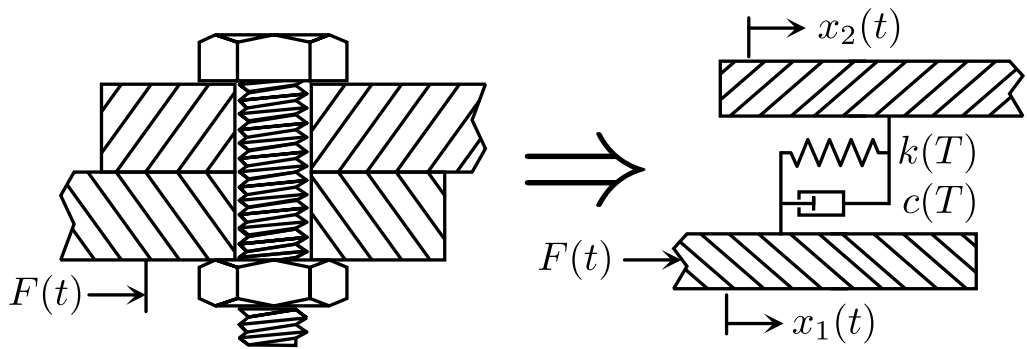

Figure 1. Schematic representation of a bolted lap joint and the proposed modeling approach.

Previous experiments [16,17] and our own experiments reported in this work demonstrated that both net loosening and net tightening can be achieved. The focus of the present research is on modeling and reproducing net loosening. In these cases, there is often an oscillation between tightening and loosening with the bolt loosening more each cycle than it retightens, leading to the observed net loss of preload. Our modeling approach hypothesizes that the instantaneous bolt tension can be decomposed into non-oscillatory and oscillatory components, which are added together to obtain the instantaneous tension. The non-oscillatory component captures the net loss of tension and is equivalent to the moving mean of the instantaneous tension. Since our focus is on modeling and predicting the net loss of tension, we focus only on the non-oscillatory component and will treat the instantaneous tension as being equivalent to this component throughout the remainder of this work.

### 2.2. Experimental System

The experimental system considered in this work is depicted in Figure 2(a) and consists of two harmonic oscillators. The first oscillator is called the large oscillator (LO) and is constructed with a stainless steel block with mass 8.515 kg of length 0.279 m, width 0.152 m, and thickness 0.0254 m. The second oscillator is called the small oscillator (SO) and is constructed with an aluminum block with mass 0.8804 kg of length 0.152 m, width 0.152 m, and thickness 0.0254 m. Each oscillator is coupled to ground (optical table) using a pair of steel flexures with width 0.152 m and length 0.114 m (active length of 0.103 m). The flexures are bolted to the optical table using aluminum L-brackets with 1/4"-20 UNC and 10-32 UNF bolts. The thickness of the flexures are 0.508 mm (0.020 in) and 1.067 mm (0.042 in) for the LO and SO, respectively, such that the flexures on the SO are significantly thicker than those for the LO. The LO is designed to have a high mass and low stiffness to ensure a low natural frequency whereas the SO is designed to exhibit a high natural frequency due to low mass and high stiffness. Both oscillators were painted white and speckled using a Sharpie® for future digital-image-correlation (DIC) testing that is not considered in this paper.

A T-shaped aluminum plate was bolted onto each oscillator using 10-32 x 0.5" UNF bolts to form the interface between the two oscillators. The T-shaped plates are shown in Fig. 2b and each have total length of 0.0826 m, a maximum width of 0.152 m, and a maximum thickness of 0.0127 m. The base of the T-shaped plates have length of 0.0762 m, width of 0.0254 m, and thickness of 0.00635 m. The interface has a length of 0.00635 m, a width of 0.0254 m. The total thickness of (two members) at the contact is 0.0127 m. The T-shaped plates were painted white and speckled using a Sharpie® for future DIC testing similar to the oscillators. A small gap between the two plates when not bolted together was incorporated to ensure that each oscillator could move without touching the other one. The T-shaped plates have mass of 0.0652 kg and 0.0654 kg respectively, bringing the total mass of the LO and SO to 8.625 kg and 0.9888 kg, respectively.

The two T-shaped plates are bolted together using a single 1/4"-20 x 1" UNC bolt (Fig. 2c) with two washers and one nut made. The bolt, washers, and nut are made from grade 8 steel with Zinc-flake-coating to ensure a material match between all three and especially for the threads. A 0.002 m diameter hole with depth of 0.012 m was drilled into the bolt head and a single HBK TB21 cylindrical half-bridge strain gauge was installed and wired following the procedure described in [63]. HBK EP70 epoxy was used to secure the strain gauge inside and was baked in an oven at 60ºC for 4 hours. After baking the epoxy, the lead and cable wires were soldered to terminals on the bolt head and a silicon gel was applied to protect them during operation. The bolts were calibrated using a tension test on an Instron 5982 universal test machine and this procedure is discussed in Section 2.4.

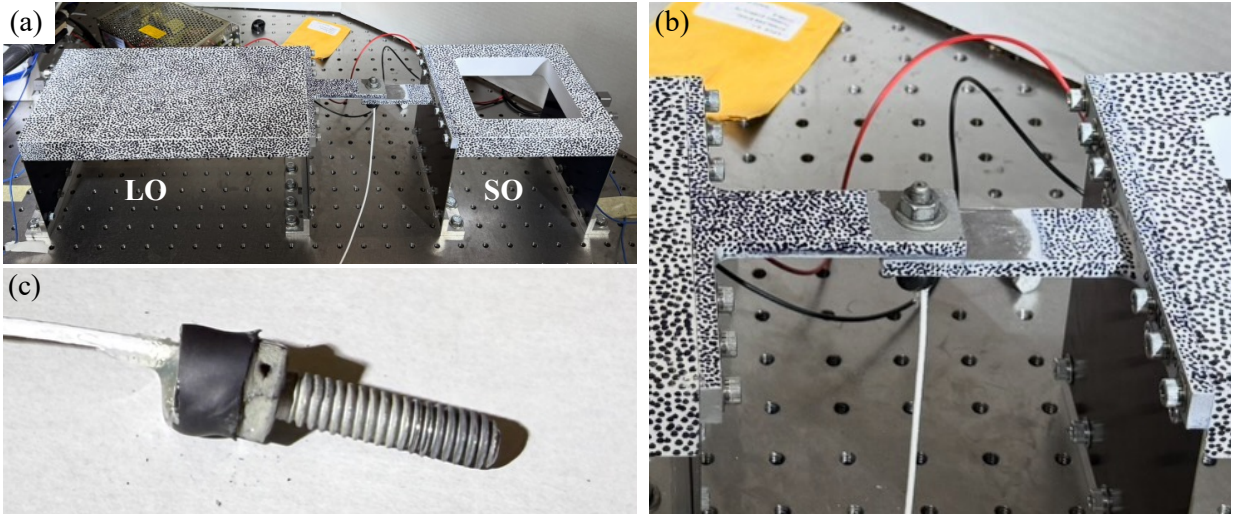

Figure 2. The experimental system. (a) View of the LO and SO coupled with the instrumented bolt. (b) Zoomed-in view showing the T-shaped plates used to create the interface. (c) The instrumented bolt.

## 2.3. Measurement Setups

All measurements focused on the transient response of the oscillators in various configurations when subjected to an impact from a modal hammer (PCB Piezotronics model 086D05) with a white, plastic tip installed. Nearly all measurements were performed with the impact applied to the LO using a custom automatic modal hammer [64,65] shown in Fig. 3. The custom automatic modal hammer uses a stepper motor (Kinco 57HBS30) controlled with the AccelStepper Library [66] and an Arduino UNO to provide impacts with repeatable amplitudes and location. A handful of measurements were performed with the impact applied manually to the SO to characterize its behavior in the uncoupled state and these are discussed in detail in Section 3.1. Each oscillator was instrumented with one high-sensitivity (PCB Piezotronics model 333B50) and one low-sensitivity (PCB Piezotronics model 353B15) accelerometer with nominal sensitivities of 102 mV/(m/s$^2$) and 1.02 mV/(m/s$^2$), respectively. The high-sensitivity accelerometers were used to measure the responses of the system when subjected to low-amplitude excitations (<100 N) whereas the low-sensitivty accelerometers were used for high-amplitude excitations (>1500 N).

The modal hammer and accelerometers were measured using a HBK QuantumX MX1601B module and the instrumented bolt was measured using a HBK QuantumX MX1616B module. The modal hammer and accelerometers were sampled at 19200 Hz (the maximum possible on the MX1601B) whereas the instrumented bolt was sampled at 4800 Hz (the maximum allowed on the MX1616B). A high-sampling rate was used for the accelerometers and modal hammer to ensure an accurate measurement of the applied force for later use in numerical simulations. The response was measured for 30.5 seconds to capture the full response of the system. A pretrigger of 0.5 seconds was used to provide enough time to measure the bolt tension prior to the impact. The measurements were captured using catman® software provided by HBK and were exported to and analyzed using MATLAB®. As needed, the accelerations were decimated in MATLAB® down to 4800 Hz to match the bolt tension measurements. The accelerations were high-pass filtered, then numerically integrated and the result high-pass filtered to obtain the velocities of the system. The velocities were then numerically integrated and high-pass filtered to obtain the displacements. The *cumtrapz* function in MATLAB® was used to perform the integration and a third-order Butterworth filter was applied with the *filtfilt* function in MATLAB® to high-pass filter the data. The *filtfilt* function is used because it applies the filter in both forwards and reverse directions, resulting in a zero-phase shift in filtered output. The cutoff frequency of the filter was 2 Hz for tests where the oscillators were not bolted together and 5 Hz for cases where they were.

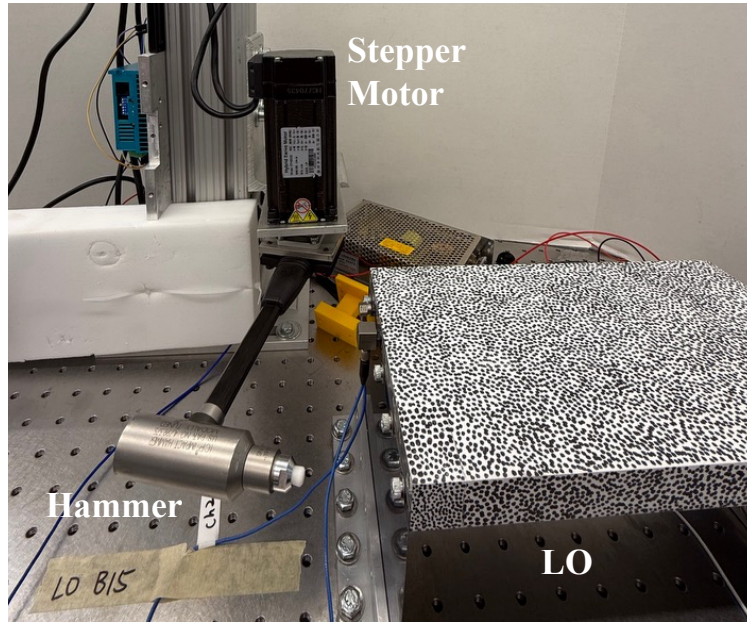

Figure 3. The automatic modal hammer used to impact the LO.

**2.4. Bolt Calibration Procedure**

The bolt was calibrated using an Instron 100 kN tensile test machine with a set of custom-made components (Fig. 4a) designed to fit the available jaws on the machine. The bottom mount is a threaded stainless steel block while the top mount is an aluminum T-bar, both speckled for DIC testing. The bolt strain was measured using the HBK data-acquisition system whereas the tensile force applied by the machine was measured using the Instron through Bluehill® Universal software. Several attempts were made to measure the applied force and resulting strains using the same system, but this was not possible due to the age and incompatibility of the Instron. Additionally, part way through the measurement, the force and measured strain would plateau and remain constant temporarily before increasing again together. The plateau was attributed to slipping of the jaws in the Instron due to the age of the system and was ignored in the calibration. The force was sampled at 100 Hz while the bolt strain was sampled at 4800 Hz, such that the strain was measured at 48 times the rate of the force.

The bolt was loaded to a maximum force of 10000 N then unloaded back to 10 N, and this cycle was repeated four times. The first cycle differed significantly from the other three, which were nearly identical in both loading and unloading responses. The difference in the first cycle was attributed to slipping of the jaws and settling of the system, and it was excluded from the calibration dataset. The remaining three cycles were averaged together and only the loading cycle was kept for the calibration. The measured force and strains were temporally aligned using the peak of the first cycle, then interpolated strain to the force time vector. The resulting force was plotted as a function of the strain and a polynomial was identified using the polyfitZero function from the MATLAB® File Exchange [67]. The resulting relationship between the force and strain is

$$T = -575.63\varepsilon^3 + 2363.7\varepsilon^2 + 2718.7\varepsilon. \tag{1}$$

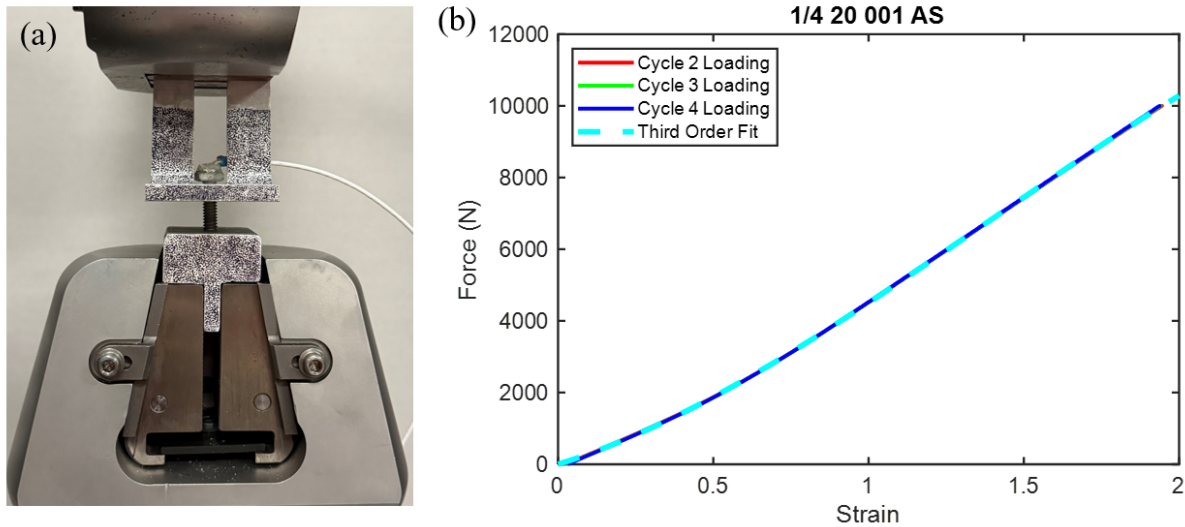

Figure 4. The calibration of the instrumented bolt. (a) The experimental setup showing the mount used to calibrate the instrumented bolt. (b) The resulting force vs strain curve and the identified polynomial model.

## 3. Experimental Identification

The experimental measurements were divided into three different types: uncoupled-oscillator tests, coupled-oscillator tests, and loosening tests. The uncoupled-oscillator tests consisted of measuring the free-response of each oscillator due to an impact from a modal hammer without the bolt installed in the interface, such that the oscillators were decoupled from each other. These measurements are used to identify the grounding stiffness and damping for each oscillator without the coupling due to the joint. The coupled-oscillator tests focused on impacting the LO with a low-amplitude force (<100 N) to induce a linear response with the bolt installed and tightened to varying tensions. The measured response is used to identify and model the effective stiffness and damping introduced by the joint for tensions varying from 6 N to 3000 N. The loosening tests consisted of impacting the LO with a large-amplitude force (>1500 N) to cause a net change in the bolt tension. These measurements were used to construct a mathematical model for the net loss of tension in the bolt due to the impact excitation. Each of these tests and the subsequent identifications are discussed in detail throughout this section.

### 3.1. Uncoupled Oscillator Tests

In the uncoupled configuration, each oscillator behaves as a damped harmonic oscillator with a single DOF as shown in the schematic representation in Fig. 5a. The experimental system in the uncoupled configuration is shown in Fig. 5b. Each oscillator was excited separately by a manual impact from a modal hammer. The automatic modal hammer was not used for these tests due to the linear behavior of each oscillator in this configuration. The resulting response was processed as discussed in Section 2.3, then the frequency-response functions (FRFs) were computed for each oscillator, which are depicted in Fig. 5c. The natural frequencies and damping ratios were identified using the rational fraction polynomial method [68], which were then used to identify the stiffness and damping coefficient for each oscillator. The natural frequencies for the LO and SO were determined to be 4.892 Hz and 44.27 Hz, respectively, and their damping ratios were determined to be 0.0155 and 0.004, respectively. These resulted in a stiffness of 8148.7 N/m for the LO and 76500 N/m for SO, and damping coefficients of 8.200 Ns/m for the LO and 2.696 Ns/m for the SO. The model FRFs are shown in Fig. 5c and show good agreement with the experimental FRFs.

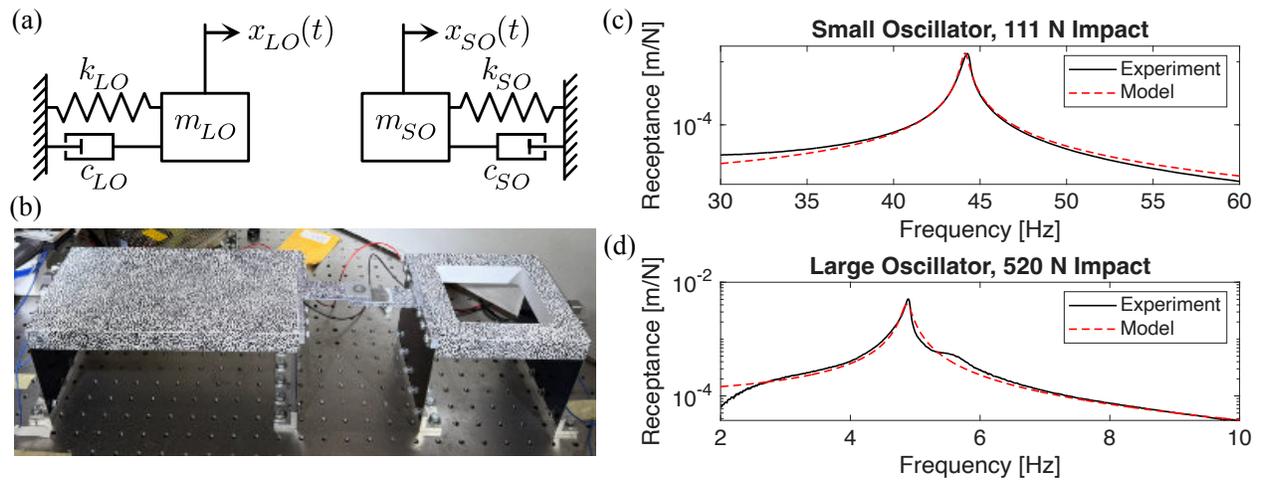

Figure 5. (a) Schematic representation of the uncoupled oscillators. (b) The experimental system in the uncoupled state. Comparisons of the measured and modeled FRFs for (c) the SO and (d) the LO.

### 3.2. Coupled Stiffness and Damping Tests & Modeling

The second set of measurements focused on measuring the linear response of the system with the bolt installed and are referred to as coupled-oscillator tests. A schematic representation of the system in this configuration is shown in Fig. 6a. In these tests, the LO was impacted using the automatic modal hammer with a low-amplitude force (<100 N) to induce a linear response. The reasons for using a low-amplitude force and inducing a linear response are to avoid altering the tension in the bolt and to avoid incurring nonlinear effects due to micro-slip or macro-slip in the joint. Consequently, the data and models produced by this study reproduce the effective linear behavior of the system and extensions to nonlinear joint models (e.g., Iwan elements [69]) are left open for future work.

An example displacement response and corresponding Fourier spectra are shown in Fig. 6b for an impact of 40 N and a preload of 1018 N, which shows that the LO and SO exhibit nearly the same motion at a single frequency around 15.4 Hz. The measured bolt tension is provided in Fig. 6c and shows that the average tension remained constant throughout the impact and motion of the oscillators. This process was repeated for 58 different initial tensions ranging from 5.8 N to 3013 N. The natural frequency of the first mode was estimated for each measurement by computing the average period between consecutive peaks in the first 3 s. Both positive and negative peaks are used to compute the average period. Due to noise and small but non-negligible ambient vibrations in the measurements, any frequencies estimated above 17 Hz were removed before averaging. The resulting natural frequencies are presented as a function of the bolt tension in Fig. 6d, which shows how the frequency rapidly increases near zero tension before plateauing around 500 N and above. Additionally, the consecutive peaks were employed with the logarithmic decrement method to compute the damping ratio of the first mode.

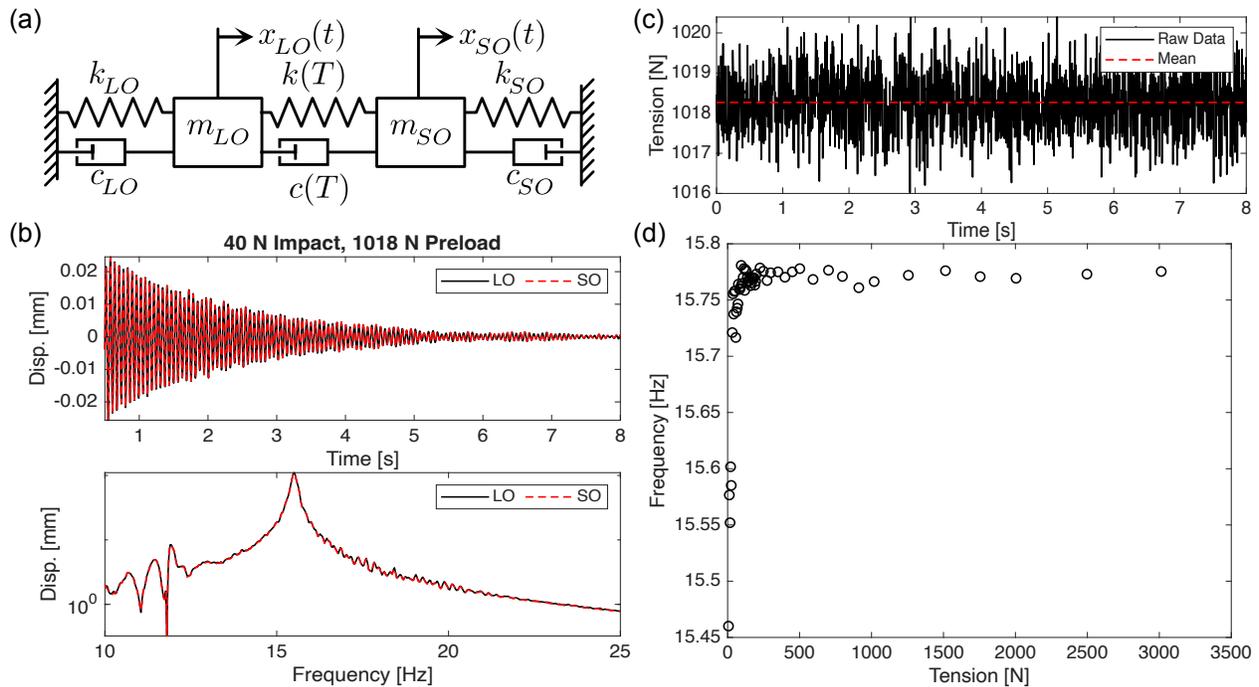

Figure 6. (a) Schematic representation of the coupled oscillators. (b) Example displacements time series and Fourier spectrum for the two oscillators under low-amplitude impacts used to identify the interface coupling stiffness and damping. (c) Example instantaneous tension showing how the tension remained constant under low-amplitude impacts (same measurement case as in b). (d) The natural frequencies estimated from the average period as a function of the bolt tension.

The reason for using the average period rather than the frequency from an FFT or FRF is because the frequency must estimated with high enough precision to capture the small changes incurred by tightening the bolt. Consequently, the response must be measured for at least 120 s to provide enough frequency resolution; however, this leads to the majority of the signal being response caused by ambient vibrations and signal noise, which leads to increased oscilations in the Fourier spectra. This issue can be overcome by applying a delayed exponential to the measured displacements without affecting the damping provided that the delay is large enough. but due to the low excitation level, this leads to increased errors and noise effects due to the majority of the signal and delayed exponential windows must be applied to enhance the clarity of the peak of the first mode. Despite these workarounds, due to the large number of measurements needed, we opted to record for a shorter duration with higher sampling frequency and determine the natural frequencies using the average periods instead.

The estimated natural frequencies are used to estimate the stiffness of the interface using the following procedure. First, the equations of motion for the coupled system are written as

$$m_{LO}\ddot{x}_{LO} + c_{LO}\dot{x}_{LO} + c_C(T)(\dot{x}_{LO} - \dot{x}_{SO}) + k_{LO}x_{LO} + k_C(T)(x_{LO} - x_{SO}) = F(t), \quad (2a)$$

$$m_{SO}\ddot{x}_{SO} + c_{SO}\dot{x}_{SO} + c_C(T)(\dot{x}_{SO} - \dot{x}_{LO}) + k_{SO}x_{SO} + k_C(T)(x_{SO} - x_{LO}) = 0, \quad (2b)$$

where $k_C(T)$ and $c_C(T)$ are the coupling stiffness and damping, respectively. From these equations, the coupling stiffness can be determined analytically by constructing the characteristic equation from the eigenvalue problem for the undamped system, then solving that for the coupling stiffness. This process results in

$$k_C = \frac{-k_{LO}k_{SO} + (k_{SO}m_{LO} + k_{LO}m_{SO})\omega_1^2 - m_{LO}m_{SO}\omega_1^4}{k_{LO} + k_{SO} - (m_{LO} + m_{SO})\omega_1^2}, \quad (3)$$

such that the coupling stiffness can be directly estimated from the measured natural frequency of the first mode. Theoretically, the maximum possible value for the first natural frequency is

$$f_1 = \frac{1}{2\pi}\sqrt{\frac{k_{LO} + k_{SO}}{m_{LO} + m_{SO}}}. \quad (4)$$

Using the values for $k_{LO}$ and $k_{SO}$ found in Section 3.1, the theoretical maximum frequency evaluates to 14.93 Hz. However, the maximum frequency estimated from the experiments is 15.78 Hz using the average period method and 15.77 Hz using the Fourier spectrum. Thus, we suspect that the flexures used to ground each oscillator experience a slight buckling as the bolt is tightened, which increases the stiffness of each oscillator. Future experiments are planned to investigate this hypothesis using digital image correlation, but are out of scope for this paper. To compensate for this, we increased the stiffness of the LO by 10% and the SO by 20%, such that $k_{LO} = 8963.6$ N/m and $k_{SO} = 91800$ N/m. This increase results in a theoretical maximum frequency of 16.29 Hz.

Using the Eq. 3 and the modified stiffnesses, we compute the coupling stiffness of the interface as a function of the tension in the bolt and depict the result in Fig. 7a. The coupling stiffnesses exhibit a similar trend to that observed for the natural frequencies in that they increase drastically until plateauing around 150 N and above. This comparable trend is expected and provides a sanity check for the results of the stiffness estimation. The coupling stiffness is modeled using the logistic function to capture the observed trend in the form

$$k_C(T) = \frac{k_I}{1 + \exp(\alpha(T - \beta))}. \tag{5}$$

To determine $k_I, \alpha,$ and $\beta$, we used the *fit* function from the *Curve Fitting Toolbox* in MATLAB® with a custom equation implementing Eq. 5, default options, and initial guesses of $k_I = 10^6$ N/m, $\alpha = -0.1$ N$^{-1}$, and $\beta = 5$ N. The *Curve Fitting T*oolbox optimizes the parameters using a nonlinear least squares algorithm and resulted in $k_I = 9.763 \times 10^5$ N/m, $\alpha = -0.0608$ N$^{-1}$, and $\beta = 2.003$ N. The identified model is shown in Fig. 7a and shows good agreement with the data with an R-squared value of 0.9355.

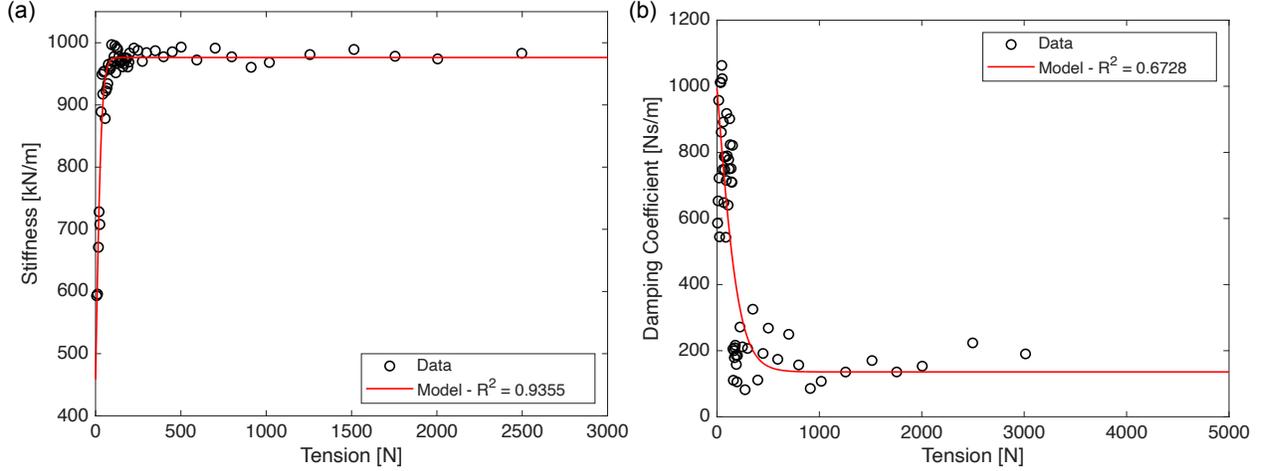

Figure 7. (a) Comparison of the interface stiffnesses computed from the experimental natural frequencies and the identified mathematical model. (b) Comparison of the experimental damping coefficients computed using logarithmic decrement and the identified mathematical model.

The logarithmic decrement method is used to estimate the damping ratio for each measurement case using the identified peaks. Since multiple peaks are available, the logarithmic decrement method is applied using the second peak as the initial point and all subsequent peaks normalized by the corresponding number of cycles. This approach provides multiple estimations of the damping ratio, which are averaged together to produce the final estimate for the response. An analytical expression relating the coupling damping coefficient to the damping ratio was constructed by solving for the first eigenvector analytically (using Mathematica), mass-orthornormalzing, and computing the modal damping coefficient of the first mode. The resulting expression for the modal damping coefficient is set equal to $2\zeta_1\omega_1$, then solved to produce an equation for the coupling damping coefficient. This approach uses the modified stiffnesses as discussed previously. Specifically, let $u_1$ be the unnormalized first eigenvector of the undamped system, then the mass-orthornormalized results in

$$\phi_1 = \frac{u_1}{\sqrt{u_1^T M u_1}}, \tag{6}$$

which is used to compute the modal damping coefficient as

$$c_1 = 2\zeta_1\omega_1 = \phi_1^T C \phi_1, \tag{7}$$

where

$$C = \begin{bmatrix} c_{LO} + c_C & -c_C \\ -c_C & c_{SO} + c_C \end{bmatrix}. \quad (8)$$

Equation 7 is solved to produce an analytical expression for the coupling damping coefficient, which results in

$$c_C = \frac{A}{B}, \quad (9)$$

where

$$A = (m_{LO}(k_C + k_{SO}) - m_{SO}(k_C + k_{LO}))(m_{LO}c_{SO} - m_{SO}c_{LO}) \\ + (4m_{LO}m_{SO}\zeta_1\omega_1 - m_{LO}c_{SO} - m_{SO}c_{LO})P, \quad (10a)$$

$$B = (m_{SO}k_{LO} - m_{LO}k_{SO})(m_{LO} - m_{SO}) - k_C(m_{LO} + m_{SO})^2 + (m_{LO} + m_{SO})P, \quad (10b)$$

$$P = \sqrt{k_C^2(m_{LO} + m_{SO})^2 + 2k_C(m_{LO} - m_{SO})(k_{SO}m_{LO} - k_{LO}m_{SO}) + (k_{SO}m_{LO} - k_{LO}m_{SO})^2}. \quad (10c)$$

The resulting coupling damping coefficients are plotted against the bolt tension in Fig. 7b, which reveals that the damping decreases significantly before plateauing after approximately 500 N. We note that the damping exhibits an opposite trend compared to the stiffness, which aligns with the expectations for a bolted joint. The damping is modeled using a reversed logistic equation in the form of

$$c_C(T) = c_D - \frac{c_I}{1 + \exp(\eta T)}. \quad (11)$$

To determine $c_D$, $c_I$, and $\eta$, we used the *Curve Fitter* toolbox in MATLAB with a custom equation implementing Eq. 11, default options, and initial gueses of $c_D = 2000$ Ns/m, $c_I = 2000$ Ns/m, and $\eta = 0.01$ N$^{-1}$. The optimization resulted in $c_D = 1853.0$ Ns/m, $c_I = 1717.4$ Ns/m, and $\eta = -0.00922$ N$^{-1}$. The identified model is shown in Fig. 7b and shows reasonable agreement with the data with an R-squared value of 0.6728. A greater variance in the damping coefficient and larger deviation in the model is expected due to significant influence of the contact region on the damping, which can vary greatly in the present system. Nevertheless, the identified mathematical model for the damping captures the overall trend of coupling damping coefficient and is suitable for use in reduced-order models. Interestingly, the model converges to a near constant value around 850 N, which is significantly higher than the convergence observed for the stiffness model.

### 3.3. Loosening Measurements and Modeling

The third set of measurements focused on inducing a change in the bolt preload by applying a relatively high-amplitude impact (average of 1608.7 N with a deviation of 81.65 N) from the automatic modal hammer to the LO. These measurements are called "loosening tests" due to the goal of causing the bolt to loosen and 84 cases were obtained in the testing. Generally, loosening tests were performed after coupled-oscillator tests, such that for each tension measured there would be both a non-loosening (linear) measurement and a loosening measurement. However, for preloads below 1100 N, the bolt tension would often increase (i.e., tighten) instead of decreasing, so loosening tests were repeated for those tensions until a decrease in preload was observed. This resulted in a cluster of measurements around certain tensions (e.g., 700 N required nine trials to observe loosening) with all but one case exhibiting tightening. The cases where tightening occurred were kept and used to demonstrate how the proposed modeling approach can also reproduce tightening if desired. However, since the focus of this work, we will later neglect the tightening cases and focuses solely loosening cases to produce a reduced-order model for the rate of tension loss. Future work is planned to investigate tightening and modify the models to capture both tightening and loosening in the relevant regimes.

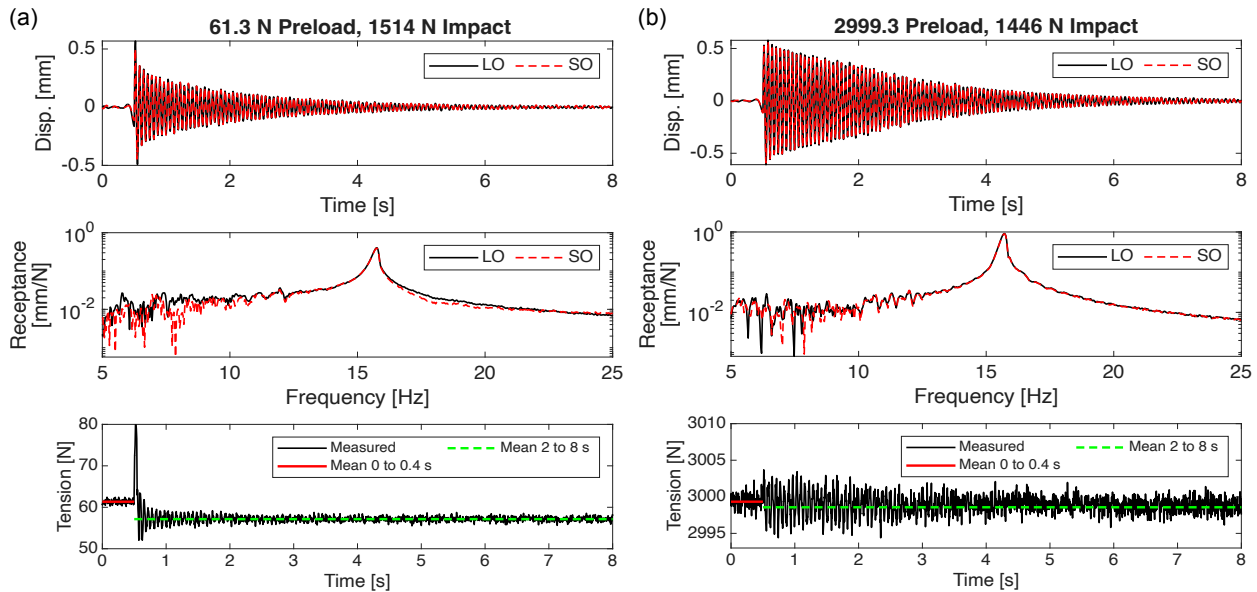

Figure 8. Example results for loosening tests for (a) low preload (61.3 N) and (b) high preload (2999.3 N). Each subfigure depicts the displacement time series and Fourier spectra for LO and SO, and the instantaneous tension with the mean values before and after the impact.

Example displacement responses, the corresponding FRFs, and the instantaneous tensions are shown in Figs. 8a and b for preloads of 61.3 N and 2999.3 N, respectively. The displacments reveal that the response of the LO and SO are nearly the same and relatively periodic despite the large-amplitude impact and the nonlinearity induced by the change in bolt tension. However, the amplitude does not exhibit a true exponential decay like the response in Fig. 6a produced in the coupled-oscillator tests. The FRFs reveals a slight asymmetry in the peak corresponding to the first mode, which implies that the response is slightly nonlinear due to the change in preload. For the instantaneous tension shown in Fig. 8a, the impact induces a sudden increase followed by a permanent decrease in the mean tension. This behavior corresponds to the non-oscillatory component described in Section 2.1. Additionally, after the sudden increase, the tension oscillates around the new mean value and decays, which corresponds to the oscillatory component discussed in Section 2.1. Considering the instantaneous tension shown in Fig. 8b, the impact results in a small change in the mean tension, while also causing the tension to oscillate.

For each measurement case, we compute the initial and final tensions by averaging the instantaneous tension for $t \in [0,0.4]$ s and $t \in [2,8]$ s, respectively. Note that the impact was applied at 0.5 s and we compute the final tension from 2 s because all non-oscillatory components settled by this time. The change in tension is computed as the difference between the initial and remaining tension and is presented in Fig. 9 where positive change indicates tightening and a negative change represents loosening. As stated previously, the bolt tension would often increase due to the impact for preloads below 1100 N, which is clearly visible in Fig. 9. The most interesting feature of Fig. 9 is the behavior is partitioned into two regimes – Regimes I and II – with the division occurring around 850 N. Below 850 N, there is a wide spread in the amount of change in tension for both loosening and tightening cases, and the amount of change increases as the initial preload decreases. This finding makes sense because at lower preloads less energy is needed to overcome sticking in the interfaces and more of the excitation energy is available to induce relative rotation between the nut and bolt. Above 850 N, the observed changes in tensions are all less than 1.25 N regardless of the initial tension as can be seen in the zoomed-in view in the top right of Fig. 9. Thus, the partition represents the initial preload where the change in tension converges to a repeatable amount. Interestingly, this convergence at 850 N is the same value that the damping model converges around even though no change in tension occurred in the measurements used to identify the damping model. This result hints at a deeper connection between the damping and change in tension, likely through the contact area and contact conditions that arise in the joint. However, it is also likely that this occurrence is a simple coincidence and further research beyond the present paper is needed to explore this connection.

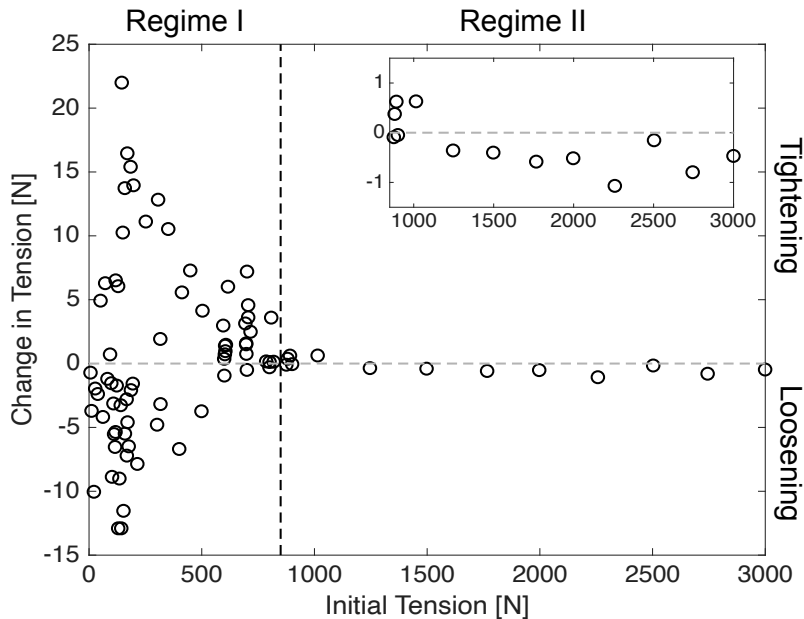

Figure 9. The measured change in tension plotted as a function of the initial tension showing that both net tightening and net loosening were achieved for preloads below 1100 N, and the behavior can be portioned into two regimes with the division at 850 N.

To illustrate the variety of behavior observed for this system, we present instantaneous tensions for a wide range of initial values in Figs. 10 and 11. Figure 10 depicts loosening cases for preloads at 100 N and below (left two columns) and for initial tensions at 1500 N and above (right two columns). Examples for preloads between these two limits are provided in Fig. 11 and are discussed in the next paragraph. The instantaneous tensions for low initial preloads exhibit a wide variety of responses whereas those for high preloads generally exhibit similar responses. The oscillatory component is generally more prevalent and observable in the low-preload cases compared to the high-preload cases. The oscillatory component is likely still present in the high preload cases but is generally smaller in amplitude than the measurement noise. The two lowest preloads exhibit strong asymmetry in their oscillations, indicating the presence of strong nonlinearity and the participation of multiple harmonics. The strong asymmetry is likely caused by impacts and pinning forces incurred by the bolt contacting the sides of the hole in the lap joints, and future research is planned to investigate and model this behavior.

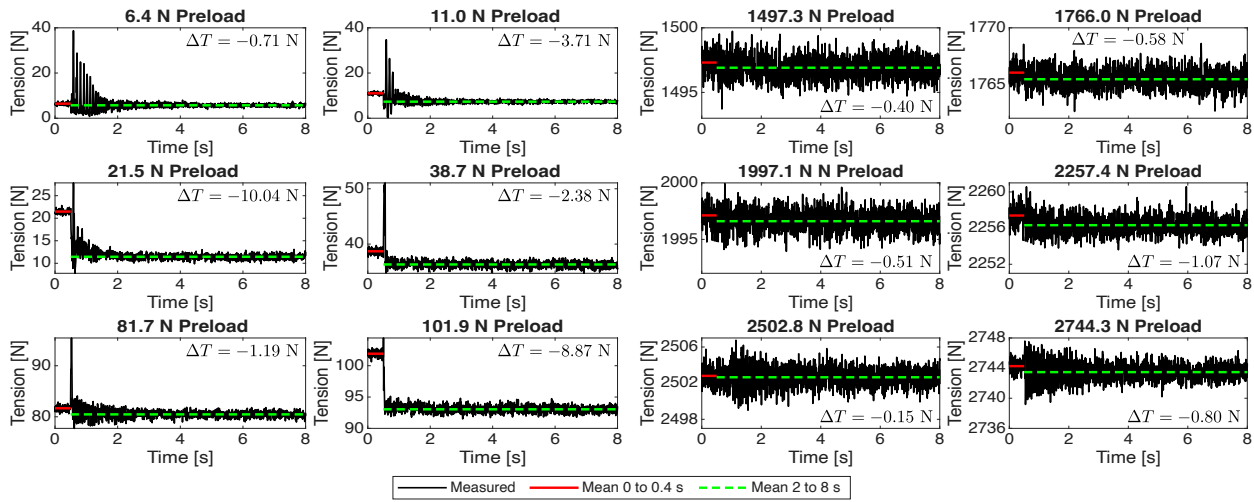

Figure 10. Example tensions where loosening occurred with the left two columns and right two columns for low and high initial tensions, respectively. The mean values before (0 to 0.4 s) and after (2 to 8 s) impact are included to show the change in tension.

Figure 11 presents example cases for preloads ranging from 120 N to 800 N. For each preload considered, we provide both net loosening and net tightening examples to show the behavior of each case. For 120 N preload (Fig. 11a), the magnitude of the change in tension is comparable for both tightening and loosening, but the oscillatory component has a significantly higher amplitude in the tightening case. Figure 11b presents the results for a preload of approximately 200 N and stark difference between the tightening and loosening cases is observed. Specifically, the change in tension is only 1.57 N in the loosening case compared to 13.96 N in the tightening case. The oscillatory component is also more prominent in the tightening case. A similar result is found for a nominal preload of 300 N (Fig. 11c), except that the oscillatory component is small in both cases. At a preload of 500 N (Fig. 11d), both cases have comparable change in tension and neither possesses a strong oscillatory component. Similar results are found for preloads of 600 N (Fig. 11e) and 800 N (Fig. 11f).

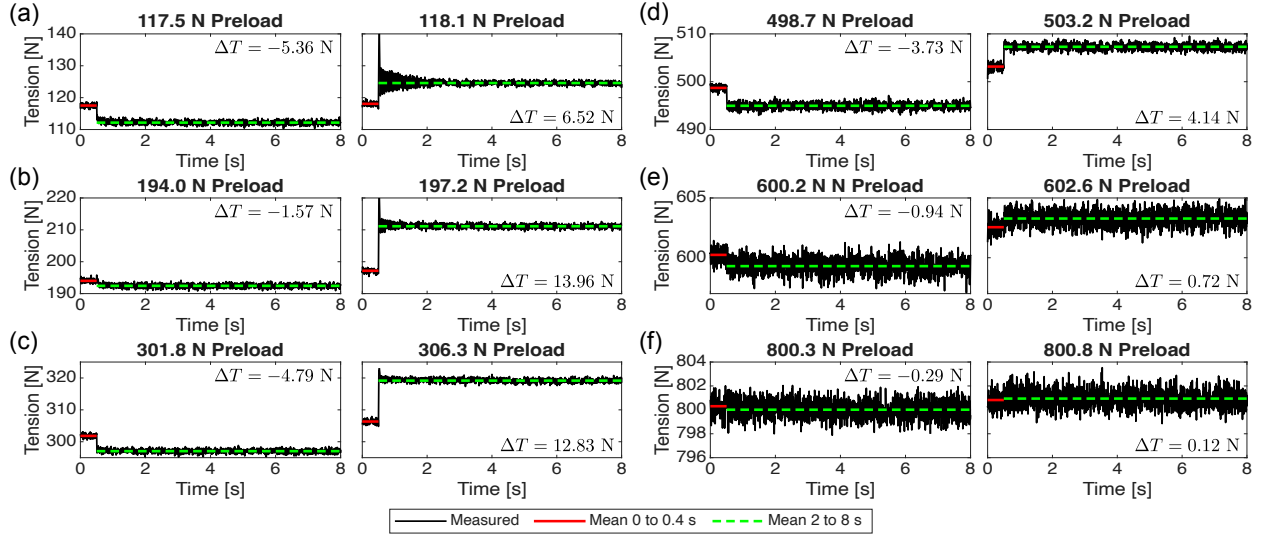

Figure 11. Example pairs of comparable tensions where loosening (left) and tightening (right) occurred for nominal preloads of (a) 120 N, (b) 200 N, (c) 300 N, (d) 500 N, (e) 600 N, and (f) 800 N. The same y-limits are used for each plot in each pair to provide a fair comparison between them. The mean values before (0 to 0.4 s) and after (2 to 8 s) impact are included to show the change in tension.

### 3.5. Loosening Model

Building on our previous work in [47,52], we model the instantaneous tension using a first-order ordinary differential equation to capture the net change in tension as

$$\dot{T} + \gamma(\dot{x}_{LO} - \dot{x}_{SO})^4 T = 0, \qquad (12)$$

where $\gamma$ is a parameter that controls the rate of change in tension. Note that the change in tension is dependent on the relative velocity and not the relative displacement. The reasoning for this is that we hypothesize that a large relative displacement can be achieved quasi-statically without inducing rotation between the bolt and the nut. Ongoing experiments are investigating this hypothesis, but those are out of the scope of the present study. We also note that the proposed model raises the relative velocity to the fourth power whereas the model used in [47,52] raised it to the second power. This change was implemented because we found that the model was able to better capture the sudden abrupt change in tension using a quartic term than with a quadratic term. Furthermore, the models used in [47,52] modeled the instantaneous torque and not the tension directly, but tension and torque are interchangeable in terms of modeling the mechanics accounting for changes in units in $\gamma$. Lastly, the instantaneous torque (or tension) was not measured directly in [47] ([52] is entirely computational) and only the initial and final torques were measured. As such, the previous model was designed to reproduce only the net change in torque and not the path or trajectory that produces that change. Since the instantaneous tensions are available in this work, we are able to adjust the proposed model to capture the observed response of the non-oscillatory component.

Our previous research [47,52] assumed that the parameter $\gamma$ was independent of tension, such that the resulting differential equation is linear. The results in [47] revealed that this is not true and this is corroborated by the trends observed in the change in tension presented in Fig. 9. However, at this stage, we treat $\gamma$ as constant and identify its value for each measurement case with negative values allowed to capture net tightening. The identification is performed by optimizing $\gamma$ until a numerical simulation of each measurement case reproduces the observed change in tension. The numerical simulation is performed with *ode45* in MATLAB® using the equations of motion in Eqs. 2a, 2b, and 12 with the tension-dependent coupling stiffness (Eq. 5) and damping (Eq. 11). The oscillators were given zero initial conditions whereas the measured starting tension was used as the initial condition for the tension model. The system was excited using the measured impact force interpolated to the integration time in the simulation. The optimization is performed with *patternsearch* in MATLAB®, which is a global optimization technique that uses a direct search method to find a true global minimum. The objective function used for the optimization is the absolute difference between the final tension measured experimentally (mean value from 2 to 8 s) and that predicted by the simulation. The mesh tolerance, function tolerance, x tolerance, maximum function evaluations, and max iterations were set to $10^{-12}$, $2.2204 \times 10^{-16}$, $2.2204 \times 10^{-16}$, $10^{10}$, and $10^{10}$, respectively. These values were chosen to force *patternsearch* to find the global minimum and are based on our previous experience with this optimization method [47,52,70,71]. The initial guess for $\gamma$ was set to $10^8$ with lower and upper bounds set to $-10^{10}$ and $10^{10}$, respectively. This range for $\gamma$ was chosen based on initial simulations with the model and made large enough to capture every possible case.

The results of the optimization are presented in Fig. 12a where positive and negative $\gamma$ values indicate net loosening and tightening, respectively. We include a zoomed-in view inlaid into the plot to show the variation in $\gamma$ at low tensions. The results show that $\gamma$ clearly varies with initial tension and we will model this behavior in the next paragraph. There is a large spread in the values for $\gamma$ below approximately 200 N with the values appearing to converge around 400 N for both loosening and tightening. However, this observation is misleading due to the large numerical values and the behavior is better depicted on a log scale. To this end, we depict the values of $\gamma$ for loosening cases only on a log y-scale in Fig. 12b, which shows that the values converge around 1100 N just like the change in tensions in Fig. 9. This outcome is not surprising given the objective function used, but provides an additional sanity check on the optimization process and methodology. Interestingly, the behavior of $\gamma$ is comparable to that observed for the coupling damping shown in Fig. 7b, except that the y-axis is on a logarithmic scale instead.

Given the similarity between the observed trends in the coupling damping and $\gamma$, we model $\gamma$ as a function of the instantaneous tension using the following equation

$$\gamma(T) = 10^{\left(\gamma_d - \frac{\gamma_I}{1+\exp(\rho T)}\right)}, \quad (13)$$

where the exponent has the same form as Eq. 11 used to model the coupling damping. The parameters were identified using the *fit* function in the *Curve Fitting* toolbox in MATLAB® with default options and initial guesses of 12, 8, and -0.0035 for $\gamma_d$, $\gamma_I$, and $\rho$, respectively. No lower or upper bounds were specified for any parameter. The initial guesses were selected from based on observations made using the custom option in the *Curve Fitting* app. The identified values are $\gamma_d = 11.79$, $\gamma_i = 7.974$, and $\rho = -0.00362$ N$^{-1}$, and the resulting model is shown in Fig. 7b. There is overall good agreement between the model and the values for $\gamma$ with the R-squared value being 0.9400.

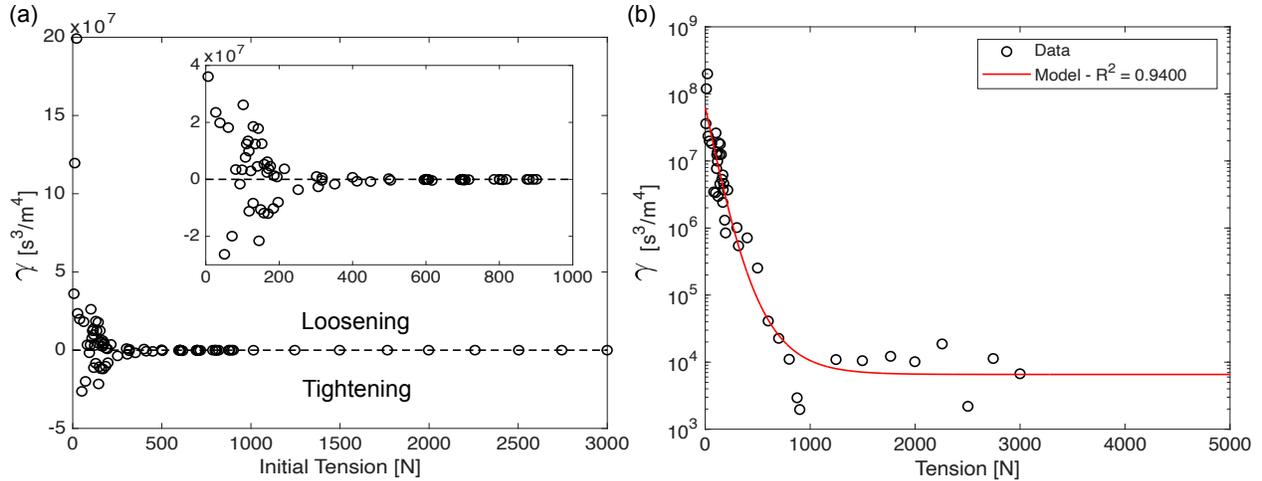

Figure 12. (a) The identified values for $\gamma$ for both tightening (negative values) and loosening (positive values) cases. (b) The mathematical model identified for $\gamma$ using only the values identified from net loosening cases.

## 4. Application of Identified Model to Reproduce Experimental Measurements

The final model for the system is formed by combining the equations of motion given by Eqs. 2a, 2b, and 12 with the models for the coupling stiffness $k_C(T)$, coupling damping $c_C(T)$, and $\gamma(T)$ provided in Eqs. 5, 11, and 13, respectively. Since the model for $\gamma(T)$ is always positive regardless of the initial tension, the final model is only able to simulate net loosening, and it fails to reproduce all cases where tightening occurs. Future research is planned to investigate what causes the system to result in a net increase in tension, but this direction is out of scope of the present work. Thus, the goal now is to determine how well the final model reproduces the experimental measurements where loosening occurred across the full range of initial tensions for both low and high excitation levels.

We start with low excitations and apply the final model to reproduce the response of the system for the coupled-oscillator tests. The system is simulated using zero initial conditions except the mean measured tension (from 0 to 8 s) is used for the instantaneous tension, and the measured impact force interpolated to the integration times is applied as the excitation to the system. We present a comparison of the measured and predicted results in Figs. 13a and b for the lowest (5.9 N) and highest preloads (3013.2 N) captured in the coupled-oscillator tests, respectively. The amplitude of the applied forces for the 5.9 N and 3013.2 N preload cases are 26.9 N and 29.1 N, respectively, such that the excitation is small enough to excite only the linear regime of the dynamics. However, the extremely small excitation means that the measured displacements are also extremely small (amplitudes of 0.02 mm) and ambient vibrations are much more apparent in the response than for higher excitations. This is evident in the experimental displacements for both preload cases through the appearance of beating patterns. Despite these effects, the model does a reasonable job of reproducing the measured displacements with the best agreement observed for the high preload. For the 5.9 N case, the simulated displacement decays faster than the experiment, which indicates that the model damping is too strong. In contrast, the simulated displacement decays at a similar rate as the experimental displacement for the 3013.2 N case, indicating that the model damping agrees with the experimental damping. This result is expected because there was a significant variation in the damping at low tension and much less variation in the damping at high tension as seen in Fig. 7b. Since the whole point of the coupled-oscillator tests was to excite the system without altering the bolt preload, the model instantaneous tension should also be constant accounting for small numerical errors. Indeed, we find that the model tension remains almost exactly the same for both cases with the difference in the initial and final model tensions are $4.74 \times 10^{-7}$ N and $3.93 \times 10^{-7}$ N for the 5.9 N and 3013.2 N cases, respectively. Although not shown here, the model was applied to all measurement cases and was able to reproduce the results without deviation. Across all measurements cases, the average difference between the initial and final tensions is $1.80 \times 10^{-6}$ N while the median is $1.58 \times 10^{-6}$ N, such that the differences are attributed to numerical error only. As such, the model is able to reproduces low-amplitude cases where no change in the bolt preload is observed or expected for the full range of preloads considered here.

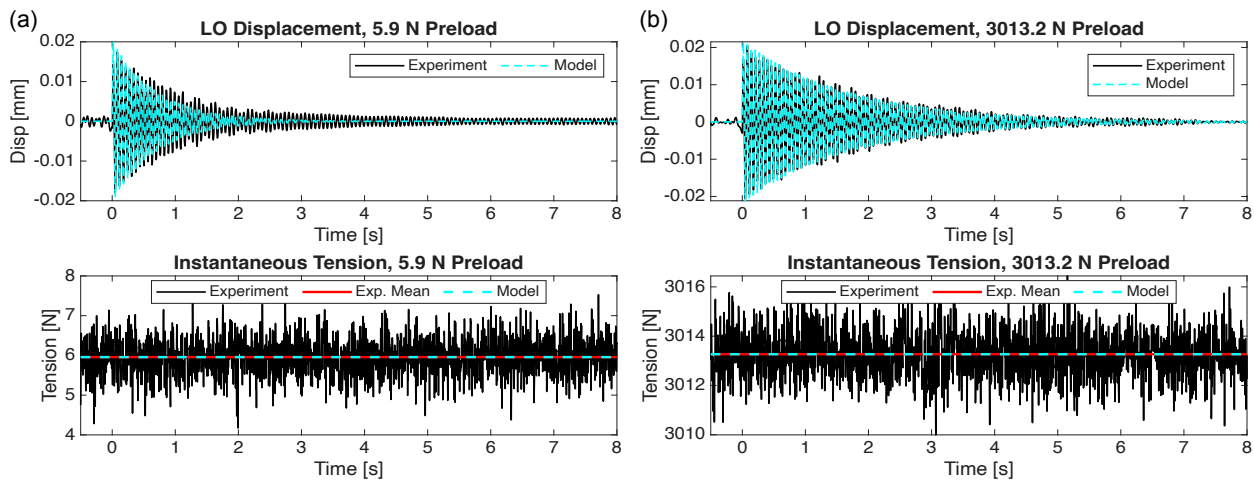

Figure 13. Comparisons between the LO displacement and the instantaneous tensions measured and predicted by the identified final model. The cases shown correspond to (a) the lowest initial tension and (b) the highest initial tension considered.

Moving on to high-excitation cases now, we apply the final model to simulate the response of the system for the loosening-test cases. We use the same type of initial conditions and excitation for the system as in low-amplitude cases, but with the initial tension and force taken from the loosening-test cases. Since the ultimate goal of the final model is to predict the loss of tension, we simulated every mesaurement case with net loosening, computed the percent error between the measured (mean from 2 to 8 s) and predicted final tensions, and present the results in Fig. 15. We find that the percent errors follow a similar trend as the changes in tension and the $\gamma$ values. Specifically, the percent errors vary significantly at low tensions and converge to a narrow band after 1100 N. This behavior makes sense and is expected because the loss of tension is primarily controlled by the value of $\gamma$, such that the region with the most variation in $\gamma$ will also produce the most variation in the percent error. The average percent error is 3.66% while the median is 1.15% The lowest percent error is $6.99 \times 10^{-5}$% with measured and predicted tensions of 2998.3 N and a difference of only -0.0021 N. The highest percent error is 62.82% with measured and predicted tensions of 11.44 N and 18.63 N, respectively, with a difference of -7.19 N.

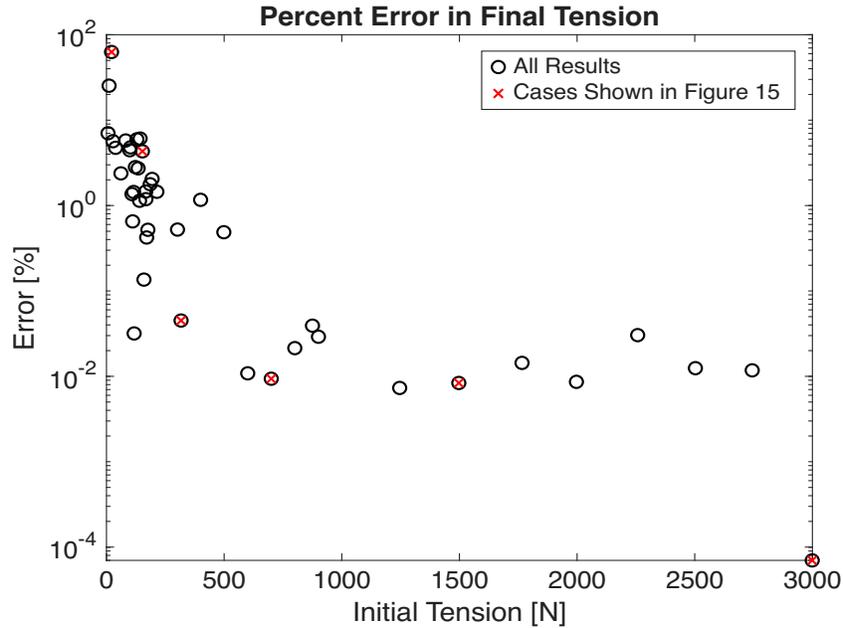

Figure 14. The percent error between the measured and predicted final tensions with specific cases shown in Fig. 15 indicated.

To illustrate the performance of the final model, we selected six cases shown with the 'x' on Fig. 15 and depict the corresponding LO displacements and instantaneous tensions in Fig. 15. We start with the lowest tension considered in the measurements in Fig. 15a with an initial preload of 6.4 N. The results show that the model does not reproduce the displacement of the LO, but does a reasonable job at reproducing the measured tension loss with a difference between the final measured and predicted tensions of 0.403 N and a percent error of 7.03%. Although not shown here, the model also does poorly at reproducing the displacement of the SO. The reason for this is that the model does not capture the strongly nonlinear dynamics that arises for low preloads and high excitation forces, such as softening due to micro-slip and pinning forces due to contact between the hole and the bolt. Consequently, the model will never be able to reproduce strongly nonlinear responses at very low preloads. Future research is planned to explore and model this regime of dynamics with a focus on how it influences the global dynamics of the parent structure.

The second case presented in Fig. 15b is for a preload of 21.5 N and this represents the case where the percent error between the measured and predicted final tensions is 62.82%, which is the highest across all cases simulated. Just like the previous case, the model is not able to reproduce the displacements of LO and SO (not shown) due to the strongly nonlinear dynamics present in this regime. Consequently, the model also underpredicts the loss in tension by -7.19 N compared to the experiment. The third case is depicted in Fig. 15c for preload of 152.7 N and the results show a good agreement for the displacement of the LO, but the model results in a loss of tension that is approximately half of that in the experiment. Indeed, the measured tension loss is 11.53 N while that predicted by the model is only 5.44 N, and the resulting percent error between the final measured and predicted tensions is 1.44%. The fourth case is for 317 N preload and is shown in Fig. 15d. This case provides an example where a relatively large loss of tension occurs and the model does well at reproducing both the displacement and instantaneous tension. The percent error in this case is only 0.045% and the predicted change in tension is 3.31 N while the measured change is 3.17 N. Figures 15e and f provide results for high preloads (700.2 N and 2999.3 N, respectively) where the model does well at reproducing both the displacements and instantaneoous tensions. Although not shown here, the model accurately reproduces both the displacements and the instantaneous tensions for all measurement cases with high preloads, as indicated by the low percent errors in Fig. 15. These results demonstrate the effectiveness of the identified final model for reproducing the experimentally measured response of the system and changes in the instantanteous tension.

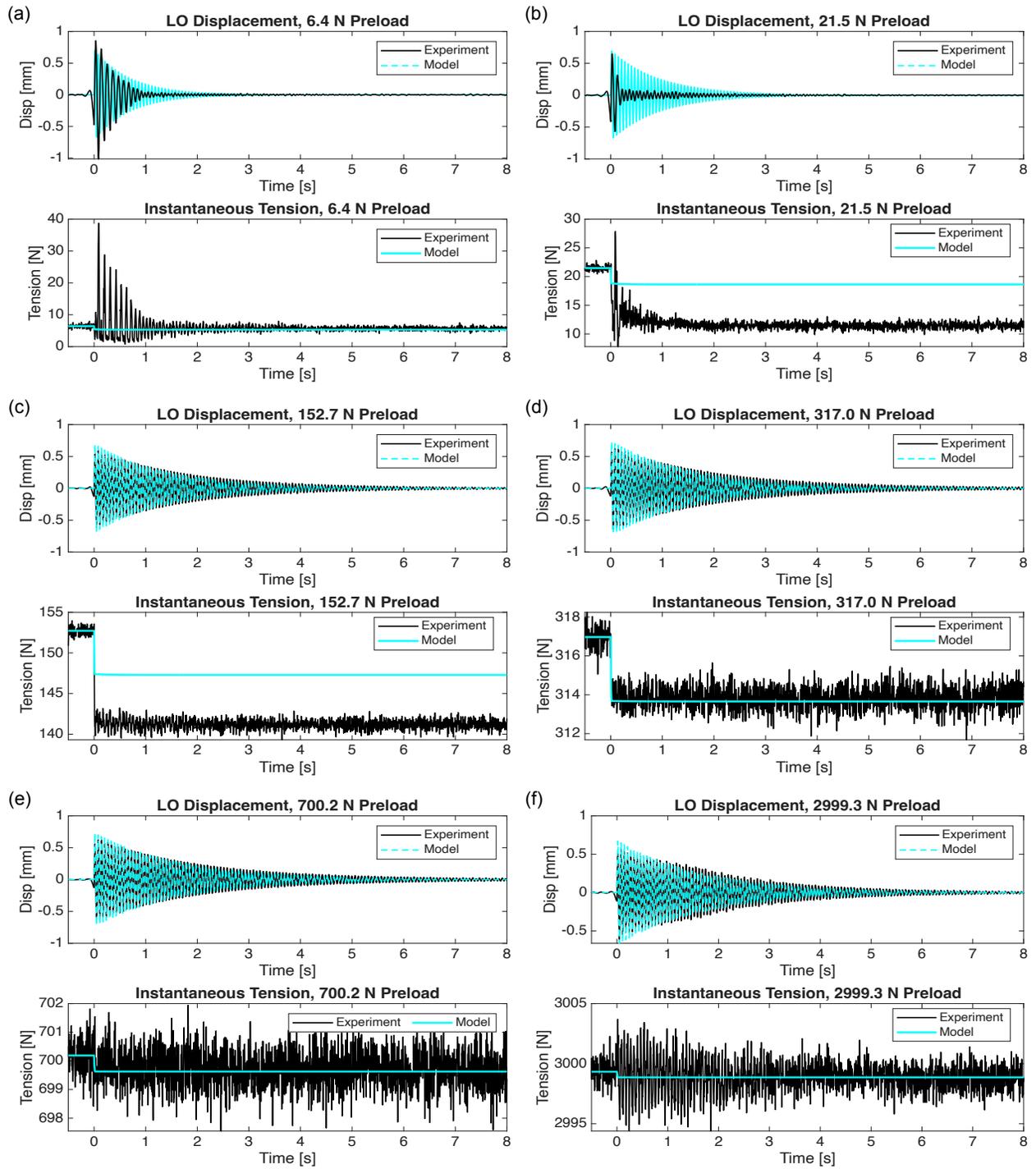

Figure 15. Comparisons between the LO displacement and the instantaneous tensions measured experimentally and predicted by the identified final model with the experimental force as input for initial preloads of (a) 6.4 N, (b) 21.5 N, (c) 152.7 N, (d) 317.0 N, (e) 700.2 N, and (f) 2999.3 N.

## 5. Concluding Remarks

This work considered the dynamics of a pair of harmonic oscillators coupled using a lap joint with a single bolt undergoing loosening. The core assertions of this work is that the bolt tension is a dynamic degree of freedom that governs the joint properties and must be modeled to capture the effects of loosening. The system was studied experimentally under three different types of measurements. First, the dynamics of each oscillator was measured and identified using transient response without the bolt installed. Second, the oscillators were coupled together with the bolt at varying preloads and excited with a low-amplitude impact from a modal hammer. The low-amplitude impact ensured that the bolt tension did not change during the motion, ensuring that the joint properties did not change either. These measurements were used to identify the linear stiffness and damping of the joint as functions of the bolt preload, and mathematical models were constructed to capture them. Third, the oscillators were excited by a large-amplitude impact in the coupled state to induce net loss of tension in the bolt, though both net tightening and loosening were observed. A first-order differential equation was selected to model the net change in tension with the rate of change of tension dependent on the relative velocity across the joint. An optimization routine was used to determine the parameter in that equation and a mathematical model was fit to capture how that parameter varies with tension. The oscillators were modeled using their equations of motion with tension-dependent stiffness and damping coupling them together. The models were combined together and used to simulate the response of the system and the instantaneous tension. The results show that the model is able to reproduce the experimental responses for a wide range of tensions provided, except for extremely low tensions and high-amplitude excitations where the behavior is strongly nonlinear. Ultimately, this research validates the critical approach of treating the bolt tension as a dynamic degree of freedom, providing a novel and effective framework for modeling and predicting the complex behavior of loosening in bolted joints.

## Funding

This research was supported by the National Science Foundation CAREER Program under grant number 2501929.

## Acknowledgements

The authors would like to express their gratitude to Aryan Singh, Felipe Kobayashi, and Blake Johnson. Their work on previous iterations of this research laid the groundwork for the success of the current project.

## Author Contributions

**Qirui He:** Conceptualization, Methodology, Software, Formal Analysis, Investigation, Data Curation, Writing – Original Draft, Writing – Review & Editing, Visualization, Supervision. **Rui Wang:** Methodology, Validation, Investigation, Data Curation, Writing – Review & Editing. **Matthew J. Alexander:** Methodology, Validation, Investigation, Data Curation, Resources,

Writing – Review & Editing. **Keegan J. Moore:** Conceptualization, Methodology, Software, Validation, Formal Analysis, Investigation, Data Curation, Writing – Original Draft, Writing – Review & Editing, Visualization, Supervision, Project Administration, Funding Acquisition